\documentclass{amsart}
\usepackage{amsmath} 
\usepackage{amssymb}
\usepackage{mathrsfs}

\newcommand{\graph}{{\rm graph}}
\newcommand{\ind}{{\rm ind}}
\newcommand{\Lim}{{\rm Lim}}
\newcommand{\Prob}{{\rm Prob}}
\newcommand{\height}{{\rm height}}
\newcommand{\mn}{{\medskip\noindent}}
\newcommand{\sn}{{\smallskip\noindent}}
\newcommand{\bbL}{{\mathbb L}}
\newcommand{\bbR}{{\mathbb R}}
\newcommand{\bbN}{{\mathbb N}}

\usepackage[hidelinks]{hyperref}

\begin{document}
\makeatletter\def\shfiuwefootnote{\gdef\@thefnmark{}\@footnotetext}\makeatother\shfiuwefootnote{Version 2015-07-03\_13. See \url{https://shelah.logic.at/papers/1061/} for possible updates.}

\title{On failure of 0-1 laws}
\author {Saharon Shelah}
\address{Einstein Institute of Mathematics\\
Edmond J. Safra Campus, Givat Ram\\
The Hebrew University of Jerusalem\\
Jerusalem, 91904, Israel\\
 and \\
 Department of Mathematics\\
 Hill Center - Busch Campus \\ 
 Rutgers, The State University of New Jersey \\
 110 Frelinghuysen Road \\
 Piscataway, NJ 08854-8019 USA}
\email{shelah@math.huji.ac.il}
\urladdr{http://shelah.logic.at}
\thanks{This work was
partially supported by European Research Council grant 338821.
\
Publication 1061 on Shelah's list.
\
The author thanks Alice Leonhardt for the beautiful typing.}





\date{July  3, 2015}

\begin{abstract}
Let $\alpha \in (0,1)_{\bbR}$ be irrational and $G_n =
G_{n,1/n^\alpha}$ be the random graph with edge probability
$1/n^\alpha$; we know that it satisfies the 0-1 law for first order logic.
We deal with the failure of the 0-1 law for 
stronger logics: $\bbL_{\infty,
{\bf k} }, {\bf k}  $ 
a large enough natural number
and the inductive logic.
\end{abstract}

\maketitle
\numberwithin{equation}{section}
\setcounter{section}{-1}

\begin{center}
\emph{Dedicated to Yuri Gurevich on the occasion of his 75th birthday}
\end{center}

\bigskip


Let $G_{n,p}$ be the random graph with set of nodes $[n] =
\{1,\dotsc,n\}$, each edge of probability $p \in [0,1]_{\bbR}$, the
edges being drawn independently, 
(see $\boxplus_1$ below). 
  On 0-1 laws (and random graphs)
see the book of Spencer \cite{Spe01} or Alon-Spencer \cite{AlSp08},
in particular on the behaviour of the random 
graph 
$G_{n,1/n^\alpha}$ for
$\alpha \in (0,1)_{\bbR}$ irrational.
On finite model theory see Flum-Ebbinghaus \cite{EbFl06}, e.g. on the
logic $\bbL_{\infty,\mathbf k}$ 
and on inductive logic,
also called LFP logic (i.e. least fix point logic). 
A characteristic example of what can 
be expressed in this logic is  ``in the 
graph 
$ G $ there is a path from the node $ x $ 
to the node $ y $", this is closed to 
what
we shall use.  
We know that $G_{n,p}$
(i.e. the case the probability $ p $ 
is constant), 
satisfies the 0-1 law for first order logic 
(proved 
independently by Fagin
\cite{Fa76} and Glebskii-et-al \cite{GKLT}).  This holds also for many
stronger logics like $\bbL_{\infty, {\bf k} } $ 
and the inductive logic.  If
$\alpha \in (0,1)_{\bbR}$ is irrational, the 0-1 law holds for
$G_{n,(1/n^\alpha)}$ and first order logic. 

The question we address is whether this holds 
also
for stronger logics as
above.  Though
our real aim is to    
address the problem for the case of graphs, the proof seems
more transparent when we have two random graph relations 
(with appropriate probabilities; we make them
directed graphs just for 
simplicity). 
  So 
we  shall deal with two
cases A and B.  In Case A, the usual graph, we have to show that there
are (just first order) formulas $\varphi_\ell(x,y)$ for 
$\ell= 
1,2$
with some special properties,
(actually we have also $ \varphi _0 ( x,y) $). 
 For Case B, those
formulas are $R_\ell(x,y),\ell=1,2$, the two directed graph relations.  
Note that 
(for Case B), the
satisfaction of the cases of the $R_\ell$ are decided directly by the
drawing and so are independent, whereas for Case A there are
(small) 
dependencies for different pairs, so the probability estimates are
more complicated.

Recall

\begin{enumerate}
\item[$\boxplus_1$]  a 0-1 context consists of:
\sn
\begin{enumerate}
\item[$(a)$]  a vocabulary $\tau$, here just the one of graphs or
  double directed graphs, 
\sn
\item[$(b)$]  for each $n,K_n$ is a set of $\tau$-models with set of
  elements = nods 
 $[n]$, in our case graphs 
or double directed graphs,
\sn
\item[$(c)$]  a distribution $\mu_n$ on $K_n$, i.e. $\mu_n:K
  \rightarrow [0,1]_{\bbR}$ satisfying $\Sigma\{\mu_n(G):G \in
  K_n\}=1$
\sn
\item[$(d)$]  the random structure is called $G_n = G_{\mu_n}$ and we tend
to speak on $G_{\mu_n}$ 
or $ G_n $ 
rather than on the context.
\end{enumerate}
\end{enumerate}
\mn 
Note that in this work ``for every random 
enough $ G_n $ \dots" is a central notion,
where:  
\mn 
\begin{enumerate}
\item[$\boxplus_2$]  for a given 0-1 context, let ``for every random
  enough $G_n$ we have $G_n \models \psi$, i.e. $G$ satisfies $\psi$"
\underline{means}  
that 
the sequence $\langle \Prob(G_n \models \psi):n \in
  \bbN\rangle$ converge to 1; of course, $\Prob(G_n \models \psi) =
  \Sigma\{\mu_n(G):G \in K_n$ and $G \models \psi\}$.
\end{enumerate}
\mn
But
\mn
\begin{enumerate}
\item[$\boxplus_3$]  $G_{n,p}$ is the case $K_n = \graph$ on 
$ [n] $ 
and
  we draw the edges independently, 
\begin{enumerate}
\item with probability $p$ when $p$
  is constant, e.g. $\frac 12$, and 
\item with probability $p(n)$ or 
probability $p_n$ when
  $p$ is a function from $\bbN$ to $[0,1]_{\bbR}$.
\end{enumerate}
\end{enumerate}
\mn
In 
the constant $ p $ 
 case, 
 the 0-1 law is strong: it is done by proving elimination
of quantifiers 
and it works also for stronger
logics:  
$\bbL_{\infty,\mathbf k}$ 
and so also for inductive logic
$\bbL_{\ind}$.  Another worthwhile case is:
\mn
\begin{enumerate}
\item[$\boxplus_4$]  $G_{n,1/n^\alpha}$ where $\alpha \in (0,1)_{\bbR}$;
  so $p_n = 1/n^\alpha$.
\end{enumerate}
\mn
Again the edges are drawn independently but the probability depends on
$n$.

The 0-1 law holds if $\alpha$ is irrational, but we have elimination
of quantifiers only up to 
(Boolean combination of)
existential formulas.  Do we have 0-1 law
also for those stronger logics? 
  We shall 
\underline{show that not}
by proving that for some so called
scheme  $ \bar{ \varphi } $ of interpretation,
for any random enough  $  G_n $,
$ \bar{ \varphi } $  
interpret an initial segment of number theory,
say up to $ m(G_n) $ 
where $ m(G_n)$ is not too small;
e.g. at least $ \log_2 (\log_2 (n)) $.

For the probabilistic argument we use 
estimates; they are as in the first order
case (see \cite{AlSp08}, so we do not repeat them).

For the full version see the author website
or the mathematical arXive.
  The 
statements for which   
we need
more estimates will probably 
 be further delayed;
those are the ones proving that:
\mn
\begin{enumerate}
\item[$\boxplus_5$]
\begin{enumerate}
\item[$\bullet$]  using $n^\varepsilon$ instead 
of $\log_2(\log_2(n))$
  in the proof 
for Case 1 so the 
value of 
``$\Prob(G_{n,1/n^\alpha}) =
  \psi$" may change more quickly,
\sn
\item[$\bullet$]  we can define ``$n$ even"
(i.e. $\Lim(\Prob(G_{n,1/n^\alpha} \models \psi$ iff $n$ is even)
  exists and is 
one;
this is done by defining a linear order on 
$ G_{ n, \bar{ \alpha } } $.
\sn
\item[$\bullet$]  we may formalize the quantification on paths, so
  getting a weak logic failing the 0-1 law, but its naturality is not
  so clear.
\end{enumerate}
\end{enumerate}
\mn
A somewhat related problem asks whether for some logic the 0-1 law
holds for $G_{n,p}$ 
(for constant $ p \in ( 0,1 )_ {\mathbb R} $, 
e.g. $p=\frac 12$) \underline{but} does not have
the elimination of quantifier, see \cite{Sh:F1166}.
\bigskip

We now try to  \underline{informally} describe the proof, 
 naturally concentrating on case B. 

Fix reals $\alpha_1 < \alpha_2$ from $(0,\frac 14)_{\bbR}$,
so $\bar\alpha = (\alpha_1,\alpha_2)$
 letting $\alpha(\ell) = \alpha_\ell$; 

\begin{enumerate}
\item[$\boxplus_6$] 
let the random 
digraph  	
 $G_{n,\bar\alpha} = ([n],R_1,R_2) = ([n],
R^{G_{n,\bar\alpha}}_1,R^{G_{n,\bar\alpha}}_2)$ with $R_1,R_2$
irreflexive relations drawn as follows:
\mn
\begin{enumerate}
\item[$(a)$]  for each $a \ne b$, we draw a truth value for
  $R_2(a,b)$ with probability $\frac{1}{n^{1-\alpha_2}}$ for yes
\sn
\item[$(b)$]  for each $a \ne b$, we draw a truth value for
$R_1(a,b)$ with probability $\frac{1}{n^{1+\alpha_1}}$ for yes
\sn 
\item[$(c)$]  those drawings are independent.
\end{enumerate}
\end{enumerate}
\mn

Now for random enough digraph 
$G= G_n= G_{n, \bar{ \alpha } }  = ([n], R_1, R_2) $ and 
node $ a \in G $ we try to  define the set
$ S_k = S_ {G,a,k} $  of nodes of $ G $
not from $ \cup \{ S_m : m< k \} $ 
by induction on $ k $ as
follows:

For $ k=0 $ let $ S_k = \{ a \} $.
Assume  $ S_0, \dots , S_k $ has been 
chosen, and we shall choose $ S_{k+1} $.

\mn
\begin{enumerate}
\item[$\boxplus_7$]
For $ \iota =1,2 $ 
 we ask: is there an $ R_\iota  $-edge $ (a,b ) $ 
with
 $ a \in S_k $
and   $ b $ not from
$ \cup \{ S_m : m \le  k \}  $?
\end{enumerate}

If the answer is no for both $ \iota =1,2 $
we stop and let ${\height}(a, G)= k $.
If the answer is yes for $ \iota= 1 $, we
 let  $  S_{k+1} $
be the set of $ b $ such that for some $ a $
the pair $ ( a,b ) $ is as above for $ \iota =1 $.,
If the answer is no for $ \iota =1 $ but yes
for $ \iota=2 $ we define $ S_{k+1} $ similarly
using  
$ \iota=2 $.

Let the height of $ G $ be 
$ \max \{ \height (a, G): a \in G \} $.
Now we can prove that for every random enough 
$ G_n $, for $ a \in G_n $ or easier- for most
$ a \in G_n $, for not too large $ k $  we have:

\mn
\begin{enumerate}
\item[$\boxplus_8$]
$ S_{G_n, a, k } $ is on the one hand  not empty
and on the other hand with $ \le n^{ 2  \alpha _2 } $
members.
\end{enumerate}

This is 
proved 
by drawing the edges not all at once
but in  $ k $ stages. In stage $ m \le k $ 
we already can compute $ S_{G_n, a, 0 } , \dots
S_{G_n, a, m } $ and 
we
have already drawn	
all the $ R_1 $-edges
and $ R_2 $-edges having at least one node
in $ S_{G_n, a, 0 } \cup  \dots \cup 
S_{G_n, a, m-1 } $;
that is for every such pair  
$ ( a, b ) $  we draw the  
truth values of  
$ R_1 (a,b ), R_2 (a,b) $. 
So arriving to $ m $ we can 
draw the edges having a nod in $ S_m $
and not dealt with earlier, and 
hence 
can compute 
$ S_{m+1} $.

The point is that in the question 
$\boxplus_7$ 
above,
if 
the answer is yes for 
$ \iota = 1 $ 
then
the number of nodes in $ {S_{m+1}} $
will be small, essentially smaller than in $
S_m $. Further,  if 
the answer for $ \iota =1 $ the answer is no
but for $ \iota =2$  the answer is yes
then 
necessarily $ S_m $ is smaller
than say $ n^{(\alpha _1+\alpha _2)/2 }  $ 
but it is known that the 
$ R_2 $-valency of of any nod of $ G_n $ is
near $ n^{ \alpha_2 } $.
So the desired inequality holds. 

By a similar argument, if we stop at $ k $
then  
in $ S_0 \cup \dots \cup S_k $ 
there are 
many nodes-
e.g.
at least near $ n^ { \alpha_2 } $
by a crud argument.
As each $ S_m $ is not too large    
necessarily   
the height of $ G_n $ is large.

The next step is to express in our logic 
the relation 
$ \{ (a_1,b_1,a_2, b_2 ) :$ for some
$ k_1, k_2 $   we have 
$ b_1  \in S_{G_n, a_1, k_1 },
b_2 \in S_{G_n, a_2, k_2 },
k_1 \le k_2 \} $.

By this we can interpret a linear order
with $ \height(G_n ) $ members.
Again using the relevant logic this suffice 
to interpret number theory up to this
height. Working more we can define a linear order 
with $ n $ elements, so can essentially
find a formula
``saying" $ n $ is even (or odd). 

For random graphs we have to work harder: instead
of having two relations we have two formulas;
one of the complications is that
their satisfaction for the relevant pairs 
are not fully independent.


\end{document}